\documentclass{amsart}

 \input xy
 \xyoption{all}

\theoremstyle{thmstyleone}%
\newtheorem{thm}{Theorem}
\newtheorem{lem}[thm]{Lemma}
\newtheorem{cor}[thm]{Corollary}
\newtheorem{prop}[thm]{Proposition}

\theoremstyle{thmstyletwo}%

\def\Cok{\mbox{\rm Cok\,}}
\def\Hom{\mbox{\rm Hom}}
\def\Ker{\mbox{\rm Ker\,}}
\def\Im{\mbox{\rm Im\,}}
\def\incl{\mbox{\rm incl}}

\begin{document}
\def\ver{s-snake-c}
\thispagestyle{empty}

\phantom m\vspace{-2cm}

\noindent{\footnotesize[{\tt \ver,} \today]}

\bigskip\bigskip
\begin{center}
  {\large\bf Snake Lemma Variations}
\end{center}

\smallskip

\begin{center}
Markus Schmidmeier

\bigskip \parbox{10cm}{\footnotesize{\bf Abstract:}
  The Kernel Complex Lemma states
  that given commutative diagram
  with exact rows and exact columns which covers the region
  under a $\Gamma$-shape,
  then the kernel sequence on the top and the kernel
  sequence at the left have in each position isomorphic
  homology.  The dual version, the Cokernel Complex Lemma,
  has been used to determine the support of a finitely presented
  functor on the Auslander-Reiten quiver.
  We note that the Snake Lemma is a consequence of the two results
  combined.
}

\medskip \parbox{10cm}{\footnotesize{\bf MSC 2010:}
  13D02, %
  16E05, 
  18G10, 
  16G70 
}

\medskip \parbox{10cm}{\footnotesize{\bf Key words:}
  Snake Lemma, Kernel Complex Lemma, Cokernel Complex Lemma
}

\bigskip

  {\it Dedicated to Professor Pu Zhang on the occasion of his 60th birthday}
\end{center}

\section{The Kernel Complex Lemma}
In this section, we show and discuss the following result.

\begin{prop}[The Kernel Complex Lemma]
  \label{prop-kcl}
  Given a commutative diagram with exact rows and columns, 
  $$\xymatrix{A_1 \ar[r]^{f_1} \ar[d]^{\alpha_1} & B_1 \ar[r]^{g_1}\ar[d]^{\beta_1}
  & C_1 \ar[r]^{h_1}\ar[d]^{\gamma_1} &\cdots\\
  A_2 \ar[r]^{f_2} \ar[d]^{\alpha_2} & B_2 \ar[r]^{g_2}\ar[d]^{\beta_2}
  & C_2 \ar[r]^{h_2}\ar[d]^{\gamma_2} &\cdots\\
  A_3 \ar[r]^{f_3} \ar[d]^{\alpha_3} & B_3 \ar[r]^{g_3}\ar[d]^{\beta_3}
  & C_3 \ar[r]^{h_3}\ar[d]^{\gamma_3} &\cdots\\
  \vdots & \vdots & \vdots & \ddots }$$
  then the two complexes given by
  taking the kernels of the vertical maps in the top row and the
  horizonal maps in the left column, 
  $$0\to \Ker\alpha_1\to \Ker\beta_1\to \Ker\gamma_1\to \cdots$$
  and
  $$0\to \Ker f_1\to \Ker f_2\to \Ker f_3\to \cdots,$$
  have in each position isomorphic homology.
\end{prop}

In the proof we use relations.  Given modules $A$, $B$,
then a {\bf relation} on $B\times A$ is a submodule of $B\oplus A$.
For example, if $f:A\to B$ and $g:A\to C$ are homomorphisms, then
$f=\{(f(a),a):a\in A\}$ is a relation on $B\times A$,
$g^{-1}=\{(a,g(a)):a\in A\}$ on $A\times C$ and
$f\circ g^{-1}=\{(b,c):(b,a)\in f\;\text{and}\;(a,c)\in g^{-1}\;\text{for some}\;
a\in A\}$ on $B\times C$.
For a proof of the Snake Lemma which uses relations we refer to \cite{rsr}.

\smallskip
We apply repeatedly the following elementary result
\cite[2.1~Lemma]{s}:

\begin{lem}
  \label{lemma-cross}
Consider the diagram with one exact row and one exact column.
$$
\xymatrix{ & B_1 \ar[d]^{\beta_1} & \\
  A \ar[r]^f & B_2 \ar[r]^g \ar[d]^{\beta_2} & C \\
  & B_3 & }
$$
\begin{enumerate}
\item Let $u\subset C\oplus B_3$ be the relation given by $g\circ \beta_2^{-1}$.
  For $(c,b_3)\in u$ we have $b_3\in \Im \beta_2 f$ if and only if 
  $c\in\Im g\beta_1$.
\item Let $v\subset B_1\oplus A$
  be the relation given by ${\beta_1}^{-1}\circ f$.
  \begin{enumerate}
    \item   Let $a\in A$.  Then $a\in\Ker \beta_2 f$ if and only if there is 
      $b_1\in B_1$ with $(b_1,a)\in v$.
    \item  Given $b_1\in B_1$, we have $b_1\in\Ker g \beta_1$
      if and only if there is 
      $a\in A$ with $(b_1,a)\in v$. \qed
  \end{enumerate}
\end{enumerate}
\end{lem}

\begin{proof}[Proof of Proposition~\ref{prop-kcl}]
  We only show that the homology at $C_0=\Ker\gamma_1$ is isomorphic
  to the homology at $K_3=\Ker f_3$.  In Part (6) we remark how the general
  situation is to be handled.
  Consider the following
  diagram where we write
  $B_0=\Ker\beta_1$, $\beta_0=\incl:B_0\to B_1$, $C_0=\Ker\gamma_1$,
  $g_0:B_0\to C_0$ for the map induced by taking kernels,
  etc.\ 
  and on the left side for $i\in\mathbb N$,
  $K_i=\Ker f_i$, $e_i=\incl:K_i\to A_i$, $\kappa_i:K_i\to K_{i+1}$.
  The diagram is commutative and all rows and columns are exact except
  possibly the top row at $C_0$ and the left column at $K_3$.
  $$\xymatrix{ & & B_0 \ar[r]^{g_0} \ar[d]^{\beta_0} & C_0 \ar[r]^{h_0}\ar[d]^{\gamma_0}
    & D_0 \ar[d]^{\delta_0} \\
    & A_1 \ar[r]^{f_1}\ar[d]^{\alpha_1} & B_1 \ar[r]^{g_1}\ar[d]^{\beta_1}
    & C_1 \ar[r]^{h_1}\ar[d]^{\gamma_1} & D_1 \\
    K_2 \ar[r]^{e_2}\ar[d]^{\kappa_2} & A_2 \ar[r]^{f_2}\ar[d]^{\alpha_2} 
    & B_2 \ar[r]^{g_2} \ar[d]^{\beta_2} & C_2 \\
    K_3 \ar[r]^{e_3}\ar[d]^{\kappa_3} & A_3 \ar[r]^{f_3}\ar[d]^{\alpha_3} & B_3 \\
    K_4 \ar[r]^{e_4} & A_4 }
  $$
  Consider the relation $u=e_3^{-1}\alpha_2f_2^{-1}\beta_1g_1^{-1}\gamma_0
  \subset K_3\times C_0$.

  \smallskip
  (1) For $c_0\in \Ker h_0$, there is $z_3\in \Ker\kappa_3$
  such that $(z_3,c_0)\in u$:
  $c_0\in\Ker h_0$ holds if and only if $c_0\in\Ker h_1\gamma_0$ since
  $\delta_0$ is a monomorphism and the rightmost square is commutative.
  By Lemma~\ref{lemma-cross} applied to the cross centered at $C_1$,
  Part (2) (a) and (b),
  there exists $b_1\in B_1$ with
  $(b_1,c_0)\in g_1^{-1}\gamma_0$;
  moreover $b_1\in\Ker\gamma_1 g_1=\Ker g_2\beta_1$.
  Again by the lemma, now centered at $B_2$, there exists $a_2\in A_2$
  with $(a_2,b_1)\in f_2^{-1}\beta_1$, and $a_2\in\Ker\beta_2f_2=\Ker f_3\alpha_2$.
  Again by the lemma, centered at $A_3$, we have $z_3\in K_3$
  such that $(z_3,a_2)\in e_3^{-1}\alpha_2$ and $z_3\in\Ker\alpha_3e_3$.
  Thus, $z_3\in \Ker \kappa_3$ since $\alpha_3 e_3=e_4\kappa_3$ and $e_4$
  is a monomorphism.
  Finally, the pair $(z_3,c_0)\in u$ since there are $b_1\in B_1$ and $a_2\in A_2$
  such that $(b_1,c_0)\in g_1^{-1}\gamma_0$, $(a_2, b_1)\in f_2^{-1}\beta_1$ and
  $(z_3,a_2)\in e_3^{-1}\alpha_2$.

  \smallskip
  (2)  Conversely, given $z_3\in \Ker \kappa_3$, there
  is $c_0\in \Ker h_0$ such that $(z_3,c_0)\in u$.

  \smallskip
  (3)  Suppose $(z_3,c_0)\in u$.  Then $c_0\in\Im g_0$ if and only if
  $z_3\in\Im \kappa_2$:
  Since $(z_3,c_0)\in u=e_3^{-1}\alpha_2f_2^{-1}\beta_1g_1^{-1}\gamma_0$,
  the pair $(e_3(z_3),\gamma_0(c_0))$ is in the relation
  $\alpha_2f_2^{-1}\beta_1g_1^{-1}$, so there exists $b_2\in B_2$ with
  $(e_3(z_3),b_2)\in\alpha_2f_2^{-1}$ and $(b_2,\gamma_0(c_0))\in\beta_1g_1^{-1}$.
  Now, $c_0\in\Im g_0$ if and only if $\gamma_0(c_0)\in\Im g_1\beta_0$ since
  $\gamma_0$ is a monomorphism.  Using Lemma~\ref{lemma-cross} (1),
  centered at $B_1$,
  this is equivalent to $b_2\in\Im \beta_1f_1=\Im f_2\alpha_1$.
  Again by this lemma, now centered at $A_2$, this is equivalent to
  $e_3(z_3)\in\Im\alpha_2e_2=\Im e_3\kappa_2$.
  Since $e_3$ is a monomorphism, this holds if and only if $z_3\in\Im\kappa_2$.

  \smallskip
  (4) There is a map $f:\Ker h_0\to \Ker\kappa_3/\Im \kappa_2$:
  Given $c_0\in\Ker h_0$, by (1) there is $z_3\in \Ker \kappa_3$ such that
  $(z_3,c_0)\in u$.  If also $z_3'\in\Ker\kappa_3$ satisfies
  $(z_3',c_0)\in u$, then $(z_3-z_3',0)\in u$.  Since $0\in\Im g_0$,
  $z_3-z_3'\in\Im \kappa_2$ by (3). Thus, the map $f$ is defined.

  \smallskip
  (5)  The map $f$ induces an isomorphism
  $\bar f:\Ker h_0/\Im g_0\cong\Ker \kappa_3/\Im \kappa_2$:
  The map $f$ is onto by (2).
  If $(z_3,c_0)\in u$ with $z_3\in\Ker \kappa_3$ and $c_0\in\Im g_0$,
  then $z_3\in\Im \kappa_2$ by (3).  Hence $f$ induces a map $\bar f$
  which must be onto.
  If $(z_3,c_0)\in u$ with $c_0\in\Ker h_0$ and $z_3\in\Im\kappa_2$,
  then $c_0\in\Im g_0$ by (2).  Thus, $\bar f$ is an isomorphism.

  \smallskip
  (6) We revisit the proof:  In order to show that homology at the third term
  $(C_0)$ in the kernel sequence on top is isomorphic to homology at the third
  term $(K_3)$ in the kernel sequence on the left, we used that the above
  diagram has 3+4 commutative squares along an antidiagonal, that the row and
  the column is exact at $B_1$, $A_2$ and at $C_1$, $B_2$, $A_3$, and that
  $\gamma_0$, $\delta_0$ and $e_3$, $e_4$ are monomorphisms.
  In general, to show for $n\in\mathbb N$ that homology at the $n$-th term
  in the kernel sequence on top is isomorphic to homology at the $n$-th term
  in the kernel sequence on the left, one will need the
  conditions on commutativitiy and exactness along the corresponding antidiagonal
  and proceed by using Lemma~\ref{lemma-cross} repeatedly.
\end{proof}

We have seen in Part (6) of the proof that certain conditions suffice
to show the partial result that the two kernel complexes
have isomorphic homology at one given position.
For the first and second terms in the kernel complex sequence we
obtain in particular:

\begin{cor}
  \label{corollary}
  Given a commutative diagram with exact rows and columns,
  $$\xymatrix{A_1 \ar[r]^{f_1} \ar[d]^{\alpha_1} & B_1 \ar[r]^{g_1}\ar[d]^{\beta_1}
    & C_1 \ar[d]^{\gamma_1} \\
    A_2 \ar[r]^{f_2} \ar[d]^{\alpha_2} & B_2 \ar[r]^{g_2}\ar[d]^{\beta_2} & C_2 \\
    A_3\ar[r]^{f_3} & B_3
  }
  $$
  then the sequences
  \begin{align*}
    0\to & \;\Ker\alpha_1\to \Ker\beta_1\to \Ker\gamma_1, \\
    0\to & \;\Ker f_1 \to \Ker f_2\to \Ker f_3
  \end{align*}
  have isomorphic homology (in positions 1 and 2). \qed
\end{cor}

\section{The Cokernel Complex Lemma}

There is a dual version of the Kernel Complex Lemma,
the Cokernel Complex Lemma.  The result has been stated and shown
in the case of $3\times 3$-regions in \cite{s}.
We state the result in the more general form and indicate
how it is used in \cite{s}
to visualize the support of a finitely presented functor.

\begin{prop}[The Cokernel Complex Lemma]
  \label{prop-ccl}
  Given a commutative diagram with exact rows and columns,
  $$\xymatrix{ \ddots & \vdots \ar[d]^{\xi_4} & \vdots \ar[d]^{\eta_4}
    & \vdots \ar[d]^{\zeta_4} \\
    \cdots  \ar[r]^{u_3} & X_3 \ar[r]^{v_3} \ar[d]^{\xi_3}
    & Y_3 \ar[r]^{w_3} \ar[d]^{\eta_3} & Z_3 \ar[d]^{\zeta_3} \\
    \cdots  \ar[r]^{u_2} & X_2 \ar[r]^{v_2} \ar[d]^{\xi_2}
    & Y_2 \ar[r]^{w_2} \ar[d]^{\eta_2} & Z_2 \ar[d]^{\zeta_2} \\
    \cdots  \ar[r]^{u_1} & X_1 \ar[r]^{v_1} 
    & Y_1 \ar[r]^{w_1} & Z_1    
  }$$
  then the two complexes given by taking
  cokernels of the horizontal maps in the right column and the vertical
  maps in the bottom row have in each component
  isomorphic homology:
  \begin{align*} \cdots\; & \to \Cok w_3 \to \Cok w_2 \to \Cok w_1 \to 0\\
    \cdots\; & \to \Cok\xi_2 \to \Cok \eta_2\to \Cok \zeta_2 \to 0
  \end{align*}
\end{prop}

The {\it proof\/} uses Lemma~\ref{lemma-cross} and is similar to the proofs of
Proposition~\ref{prop-kcl} and \cite[Proposition~1.4]{s}
and will be omitted. \qed

\bigskip
The situation in the Cokernel Complex Lemma arises in \cite{s}
where  a short exact
sequence and a short right exact sequence are given, say
$$\mathcal A: 0\to A\to B\stackrel g\to C\to 0\qquad \text{and}\qquad
\mathcal E: X\stackrel u\to Y\to Z\to 0.$$
The covariant and the contravariant $\Hom$-functors
given by the modules in the two sequences yield the following
commutative diagram with exact rows and columns, where we write
$(X,A)$ for $\Hom(X,A)$.
$$\xymatrix{ & 0 \ar[d] & 0 \ar[d] & 0\ar[d] \\
  0\ar[r] & (Z,A) \ar[r]\ar[d] & (Z,B) \ar[r]^{(Z,g)}\ar[d] & (Z,C)\ar[d] \\
  0\ar[r] & (Y,A) \ar[r]\ar[d]^{(u,A)} & (Y,B) \ar[r]^{(Y,g)}\ar[d]^{(u,B)}
  & (Y,C)\ar[d]^{(u,C)} \\
  0\ar[r] & (X,A) \ar[r] & (X,B) \ar[r]^{(X,g)} & (X,C) }
$$
According to Proposition~\ref{prop-ccl}, the two complexes,
\begin{align*} 0\to \; & \Cok(Z,g)\to \Cok(Y,g)\to \Cok(X,g)\to 0,\\
  0\to \; & \Cok(u,A)\to \Cok(u,B)\to \Cok(u,C)\to 0
\end{align*}
have in each position isomorphic homology.

\smallskip
In \cite{s} we are interested in the situation where $\mathcal A$
is an Auslander-Reiten sequence in a category of modules over a
finite dimensional algebra, see \cite{ars} for definition and basic results.
Then the dimension for $\Cok(X,g)$
determines the multiplicity of $C$ as a direct summand of $X$.
On the other hand, if
$$0\to (Z,-)\to (Y,-)\stackrel{(u,-)}\longrightarrow (X,-)\to E\to 0$$
is the projective resolution of a finitely presented covariant functor $E$,
then for an object $A$, the dimensions for $E(A)$ and for $\Cok(u,A)$
agree.
We deduce from the interplay between the homologies
of the two cokernel complexes in particular that the functor $E$ is additive
on an Auslander-Reiten sequence $\mathcal A$ (in the sense that
$\dim E(B)=\dim E(A)+\dim E(C)$) unless $C$ occurs as a direct summand
of $X\oplus Y\oplus Z$. 
In \cite{s}
we can visualize the support of the finitely presented functor $E$
on the Auslander-Reiten quiver as a hammock.

\section{The Snake Lemma revisited}

Certainly, the Snake Lemma has been shown often enough since it appeared
in \cite[Proposition~2.10]{am}.  We provide yet another proof
as the result
follows immediately from Proposition~\ref{prop-kcl} (or its corollary)
and Proposition~\ref{prop-ccl}.

\begin{prop}[The Snake Lemma]
  Given a commutative diagram with exact rows and columns,
  $$\xymatrix{ & A\ar[r]^f \ar[d]^\alpha & B\ar[r]^g \ar[d]^\beta
    & C\ar[r] \ar[d]^\gamma & 0\\
    0\ar[r] & A'\ar[r]^{f'} & B'\ar[r]^{g'} & C' & }
  $$
  there is a homomorphism $\delta:\Ker\gamma\to\Cok\alpha$ such that the
  sequence of kernel and cokernel maps is exact.
  $$\Ker\alpha\to \Ker\beta\to \Ker\gamma\stackrel\delta\to
  \Cok\alpha\to\Cok\beta\to\Cok\gamma$$
  Moreover, $f$ is monic if and only if the map $\Ker\alpha\to\Ker\beta$ is;
  and $g'$ is onto if and only if the map $\Cok\beta\to\Cok\gamma$ is.
\end{prop}

\begin{proof}
  It follows form the assumptions that the following diagram 
  $$\xymatrix{ A \ar[r]^f \ar[d]^\alpha & B\ar[r]^g \ar[d]^\beta 
    & C \ar[d]^{\gamma} \\
    A'\ar[r]^{f'} \ar[d]^{\alpha'} & B'\ar[r]^{g'}\ar[d]^{\beta'}
    & C' \\
    \Cok\alpha \ar[r]^{f''} & \Cok\beta
    }
  $$
  is commutative and has exact rows and columns.
  We obtain from Corollary~\ref{corollary} that the kernel complexes
  \begin{align*}(*)\qquad & 0\to \Ker\alpha\stackrel{\tilde f}\to
    \Ker\beta\stackrel{\tilde g} \to \Ker\gamma, \\
  &  0\to \Ker f\stackrel{\tilde \alpha}\to \Ker f'
  \stackrel{\tilde\alpha'} \to \Ker f''
  \end{align*}
  have isomorphic homology (in positions one and two).
  In particular, if $f$ is a monomorphism then $\Ker f=0$
  and the complex $(*)$ is
  exact at $\Ker\alpha$.
  Recalling the assumption that $f'$ is a monomorphism,
  we obtain in general that $\Ker \tilde f\cong \Ker f$.
  Since $f'$ is a monomorphism, 
  we also obtain that the complex $(*)$ is exact at $\Ker\beta$.  

  \smallskip
  We have seen that the following diagram 
  $$\xymatrix{& 0\ar[d] & 0\ar[d] & 0\ar[d] & 0\ar[d]  \\
    0\ar[r] & \Ker\tilde f\ar[r]\ar[d]^{\cong} & \Ker\alpha\ar[r]^{\tilde f} \ar[d]
    & \Ker\beta\ar[r]^{\tilde g}\ar[d]
    & \Ker\gamma  \ar[d]  \\
    0\ar[r] & \Ker f\ar[r] \ar[d]& A \ar[r]^f \ar[d]^\alpha
    & B\ar[r]^g \ar[d]^\beta & C \ar[d]^\gamma \\
    0\ar[r] & 0 \ar[r] & A'\ar[r]^{f'} & B'\ar[r]^{g'} & C' }
  $$
  is commutative with exact rows and columns.
  Applying Proposition~\ref{prop-ccl}, the cokernel complexes,
  \begin{align*} (**) \qquad & 0\to \Cok\alpha\stackrel{f''}\to \Cok\beta
  \stackrel{g''}\to \Cok\gamma\to 0,\\
  & 0\to \Cok\tilde g\to \Cok g\to \Cok g'\to 0
  \end{align*}
  have isomorphic homology.  Starting from the right we read off:
  Whenever $g'$ is onto,
  then $(**)$ is exact at $\Cok\gamma$.
  In general, since $\Cok g=0$ by assumption, there is an isomorphism $\Cok g''\cong \Cok g'$.
  Again since $\Cok g=0$, we have that the complex $(**)$ 
  is exact at $\Cok\beta$.  
  Finally, homology at $\Cok\tilde g$
  (which is $\Cok\tilde g$) is isomorphic to homology at $\Cok\alpha$
  which is $\Ker f''$ --- this isomorphism yields the desired connecting
  homomorphism $\delta$.
\end{proof}

\medskip\noindent
    {\bf Conflict of Interest Statement.}
    The author declares that there is no conflict of interest.
    There are no datasets associated to this paper.
    

\bigskip\noindent
Mathematical Sciences\\ Florida Atlantic University\\
Boca Raton, Florida 33431\\ {\tt\small markusschmidmeier@gmail.com}

\end{document}